\documentclass[a4paper,12pt]{article}
\usepackage{latexsym}		
\usepackage{amsmath}
\usepackage{amsthm}

\usepackage{graphicx}
\usepackage{txfonts}
\usepackage{bm}
\usepackage{color}
\usepackage{epic,eepic}

\newcommand{\RED}[1]{{\color{red}#1}} 
 \renewcommand{\RED}[1]{{#1}} 
\newcommand{\BLU}[1]{{\color{blue}#1}} 
 \renewcommand{\BLU}[1]{{#1}} 

\usepackage{geometry}
\newcommand{\memo}[1]{{\bf \small \RED{[MEMO:}} \BLU{#1} \ {\bf \small \RED{:end]} }}  
   \renewcommand{\memo}[1]{}           
\newcommand{\OMIT}[1]{{\bf [OMIT:} #1 \ {\bf --- end OMIT] }}  
   \renewcommand{\OMIT}[1]{}            

\numberwithin{equation}{section}

\begin{document}

\title{
Double-Exponential transformation:
\\ 
A quick review of a Japanese tradition
}

\author{
Kazuo Murota%
\thanks{
The Institute of Statistical Mathematics,
Tokyo 190-8562, Japan; 
and
Faculty of Economics and Business Administration,
Tokyo Metropolitan University, 
Tokyo 192-0397, Japan,
murota@tmu.ac.jp}
\ and \
Takayasu Matsuo%
\thanks{
Department of Mathematical Informatics,
Graduate School of Information Science and Technology,
University of Tokyo, Tokyo 113-8656, Japan 
matsuo@mist.i.u-tokyo.ac.jp}
}

\date{January 5, 2023 / May 21, 2025}

\maketitle

\begin{abstract}
This paper
is a short introduction 
to numerical methods using the double exponential (DE) transformation,
such as tanh-sinh quadrature and DE-Sinc approximation.
The DE-based  methods for numerical computation 
have been developed intensively in Japan
and the objective of this paper
is to describe their history 
in addition to the underlying mathematical ideas.
\end{abstract}



\section{Introduction}

The double exponential (DE) transformation
is a generic name of 
variable transformations
(changes of variables)
used effectively in
numerical computation on analytic functions,
such as numerical quadrature and function approximation.
A typical DE transformation
is a change of variable $x$ to 
another variable $t$ 
by $x = \phi(t)$ with the function
\begin{equation} \label{tanhsinhTrans}
 \phi(t) = \tanh \left( \frac{\pi}{2} \sinh t \right) .
\end{equation}
The term ``double exponential" refers to the property 
that the derivative  
\begin{equation} 
 \phi'(t) 
 =\frac { \frac{\pi}{2} \cosh t } {\cosh^{2}(\frac{\pi}{2}\sinh t)}
\end{equation}
decays double exponentially 
\begin{equation} \label{tanhsinhDerDecay}
 \phi'(t)  \approx  \exp 
\left( -\frac{\pi}{2} \exp |t|  \right)
\end{equation}
as $|t| \to \infty$.

This paper is a short introduction 
to numerical methods using DE transformations
such as 
the double exponential formula (tanh-sinh quadrature)
for numerical integration
and the DE-Sinc method for function approximation. 
The DE-based  methods for numerical computation 
have been developed intensively in Japan
\cite{Mor05,MS01,SM04sinc,TO23bk},
and a workshop
titled ``Thirty Years of the Double Exponential Transforms"
was held at RIMS (Research Institute for Mathematical Sciences, Kyoto University) 
on September 1--3,  2004
\cite{OS05deYears}.
The objective of this paper
is to describe the history of the development
of the DE-based  methods
in addition to the underlying mathematical ideas.

Originally, this paper was written
in memory of Professors Masao Iri, Masatake Mori, and Masaaki Sugihara,
and published as an article in Tokyo Intelligencer (pp.~15--21) at the
10th International Congress on Industrial and Applied Mathematics,
2023 (ICIAM 2023, Tokyo).
The present paper is its minor revision.

\section{DE formula for numerical integration}

The DE formula for numerical integration
invented by Hidetosi Takahasi and Masatake Mori 
\cite{TM74}
was first presented at the RIMS workshop 
``Studies on Numerical Algorithms,'' 
held on October 31--November 2, 1972.
The celebrated term ``double exponential formula"
was proposed there,
as we can see in the proceedings paper 
\cite{TM73defirst}.

\subsection{Quadrature formula}

The DE formula was motivated by the fact that the trapezoidal rule
is highly effective for integrals 
over the infinite interval $(-\infty,+\infty)$.
For an integral
\begin{equation}   \label{int01fx} 
I = \int_{-1}\sp{1} f(x) {\rm d}x,
\end{equation}
for example,
we employ a change of variable
$x = \phi(t)$
using some function $\phi(t)$
satisfying 
$\phi(-\infty)= -1$ and $\phi(+\infty)= 1$,
and apply the trapezoidal rule 
to the transformed integral
\begin{equation} 
I = \int_{-\infty}\sp{+\infty} f(\phi(t)) \phi'(t){\rm d}t,
\end{equation}
to obtain an infinite sum of discretization 
\begin{equation} \label{Ihinfsum}
 I_{h} = h \sum_{k=-\infty}\sp{\infty} f(\phi (kh)) \phi'(kh).
\end{equation}
A finite-term approximation to this infinite sum
results in an integration formula
\begin{equation} \label{Ihfinsum}
 I_{h}\sp{(N)} = h \sum_{k=-N}\sp{N} f(\phi (kh)) \phi'(kh) .
\end{equation}
Such combination of the trapezoidal rule with 
a change of variables was conceived by 
several authors \cite{IMT70,Sch69,Ste73,TM73} around 1970.

The error
$I - I_{h}\sp{(N)}$
of the formula \eqref{Ihfinsum} consists of two parts, the error 
$E_{\rm D} \equiv  I  -  I_{h}$
incurred by discretization \eqref{Ihinfsum}
 and the error 
$E_{\rm T} \equiv  I_{h} -  I_{h}\sp{(N)}$
caused by truncation of an infinite sum $I_{h}$ 
to a finite sum $I_{h}\sp{(N)}$.

The major findings of Takahasi and Mori 
consisted of two ingredients.
The first was that the double exponential decay 
of the transformed integrand
$f(\phi (t)) \phi'(t)$
achieves the optimal balance (or trade-off) between
the discretization error $E_{\rm D}$ 
and the truncation error $E_{\rm T}$.
The second finding was that a concrete choice of
\begin{equation} \label{tanhsinhTrans2}
 \phi(t) = \tanh \left( \frac{\pi}{2} \sinh t \right)
\end{equation}
is suitable for this purpose
thanks to the double exponential decay
shown in \eqref{tanhsinhDerDecay}.
With this particular function $\phi(t)$
the formula \eqref{Ihfinsum} reads
\begin{align}
 I_{h}\sp{(N)} = h  \sum_{k=-N}\sp{N} 
f\left(\tanh \left( \frac{\pi}{2} \sinh (kh) \right)\right)
\frac { ({\pi}/{2}) \cosh (kh) } {\cosh^{2}(({\pi}/{2})\sinh (kh))} ,
\end{align}
which is sometimes called ``tanh-sinh quadrature."
The error of this formula is estimated roughly as 
\begin{equation} \label{errorDEformula}
\big| I - I_{h}\sp{(N)} \big| \approx \exp ( -C N/ \log N )
\end{equation}
with some $C > 0$.
The DE formula has an additional feature that
it is robust against end-point singularities of integrands.

The idea of the DE formula can be 
applied to integrals
over other types of intervals of integration.
An integral 
$I = \int_{a}\sp{b} f(x) {\rm d}x$ 
over a general finite interval $[a,b]$
can be reduced to the form of \eqref{int01fx}
by a linear tranformation 
$x=[(b-a) \hat{x} + (b+a)]/2$
of the variable.
For integrals over infinite intervals, we use the following transformations:
\begin{align}
&I = \int_{0}\sp{+\infty} f(x) {\rm d}x,
\qquad
x = \exp \left( \frac{\pi}{2} \sinh t \right),
\label{DEformula0Inf}
\\
&I = \int_{-\infty}\sp{+\infty} f(x) {\rm d}x,
\qquad
x = \sinh \left( \frac{\pi}{2} \sinh t \right).
\label{DEformulaInfInf}
\end{align}
Such formulas are also referred to as the double exponential formula.
The DE formula is available in 
Mathematica (NIntegrate),
Python library SymPy,
Python library mpmath,
C++ library Boost,
Haskell package integration, etc.

\subsection{Optimality}

\begin{figure}
\centering
\includegraphics[width=0.7\textwidth]{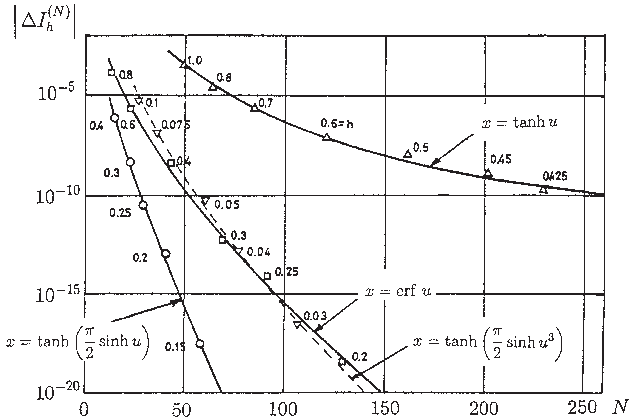}
\caption{Comparison of the efficiency of several variable transformations for
the integral
$I = \int_{-1}\sp{1} {\rm d}x/\{ (x-2)(1-x)\sp{1/4}(1+x)\sp{3/4} \}$.
The figure is taken from Mori \cite[Fig.~4]{Mor05}
with permission from 
Publications RIMS, Research Institute for Mathematical Sciences, Kyoto University;
$u$ and $N$ in the figure correspond,
respectively, to $t$ and $2N+1$ in the present notation.}
\label{FGdeintMori05}
\end{figure}

Optimality of the DE transformation 
\eqref{tanhsinhTrans2}
was discussed already 
by Takahasi and Mori \cite{TM74}.
Numerical examples also support its optimality.
Figure \ref{FGdeintMori05} 
(taken from \cite{Mor05}) shows
the comparison of the DE transformation
\eqref{tanhsinhTrans2}
against other transformations
\begin{align*}
 \phi(t) &= \tanh t ,
\\
 \phi(t) &= \tanh \left( \frac{\pi}{2} \sinh t\sp{3} \right),
\\
 \phi(t) &= {\rm erf} (t) 
  = \frac{2}{\sqrt{\pi}}
 \int_{0}\sp{t} \exp ( - s\sp{2} ) {\rm d}s
\end{align*}
for 
$I = \int_{-1}\sp{1} {\rm d}x/\{ (x-2)(1-x)\sp{1/4}(1+x)\sp{3/4} \}$
that has integrable singularities
at both ends of the interval of integration.
The DE formula converges much faster than others.
It is known that the tanh-rule 
(using $\phi(t) = \tanh t$)
has the (rough) convergence rate  
$\exp ( -C \sqrt{N} )$,
in contrast to 
$\exp ( -C N/ \log N )$ 
in \eqref{errorDEformula} of the DE formula.

The optimality argument of \cite{TM74},
based on complex function theory,
was convincing enough for the majority of scientists 
and engineers, but not perfectly satisfactory
for theoreticians.
Rigorous mathematical argument for optimality of the DE formula
was addressed by Masaaki Sugihara 
\cite{Sug86,Sug88,Sug97optDE}
in the 1980s and 1990s 
in a manner comparable to Stenger's framework \cite{Ste78}
for optimality of the tanh rule.
It is shown in \cite{Sug97optDE} (also \cite{TSMM09fnclass})
that the DE formula is optimal
with respect to a certain class (Hardy space) of integrand functions.

In principle, for each class of integrand functions we may be able to 
find an optimal quadrature formula,
and the optimal formula naturally depends on our choice of the admissible
class of integrands.
Thus the optimality of a quadrature formula is only relative. 
However, 
it was shown by Sugihara 
that no nontrivial class of integrand functions
exists that admits a quadrature formula with smaller errors than the DE formula.
We can interpret this fact as the absolute optimality of the DE formula.

\subsection{Fourier-type integrals}

For Fourier-type integrals such as
\begin{equation*}  
 I = \int_{0}\sp{+\infty} f_{1}(x) \sin x  \ {\rm d}x ,
\end{equation*}
the DE formula like \eqref{DEformula0Inf} 
is not very successful.
To cope with Fourier-type integrals,
a novel technique, in the spirit of DE transformation,
was proposed by Ooura and Mori \cite{OM91osc,OM99four}.   
In \cite{OM91osc} they proposed to use
\[
 \phi(t) = \frac{t}{ 1 - \exp(- K \sinh t) }
\]
($K > 0$),
which maps $(-\infty,+\infty)$ to $(0,+\infty)$
in such a way that
(i) $\phi'(t) \to 0$ double exponentially as $t \to -\infty$
and 
(ii) $\phi(t) \to t$ double exponentially as $t \to +\infty$.
The proposed formula changes the variable by
$x = M \phi(t)$
to obtain
\begin{align*} 
I =  M \int_{-\infty}\sp{+\infty} \!\! f_{1}(M \phi (t))
 \sin(M \phi(t)) \phi'(t)  {\rm d}t,
\end{align*}
to which the trapezoidal rule with equal mesh $h$ is applied,
where $M$ and $h$ are chosen to satisfy $Mh = \pi$.
The transformed integrand decays double exponentially
toward $t \to -\infty$
because of the factor $\phi'(t)$
and also 
toward $t \to +\infty$
because 
$M \phi(t)$ 
for $t=kh$ (sample point of the trapezoidal rule) tends
double exponentially 
to $Mt = M kh = k \pi$, at which sine function vanishes.
Another (improved) transformation function
\[
 \phi(t) 
= \frac{t}{ 1 - \exp( -2t - \alpha (1 - {\rm e}\sp{-t})  - \beta ({\rm e}\sp{t} - 1) )},
\]
is given in \cite{OM99four}, where
$\beta = 1/4$ and
$\alpha = \beta/ \sqrt{1  +  M \log (1+M)/(4 \pi)}$.

\subsection{IMT rule}

In 1969, prior to the DE formula,
a remarkable quadrature formula was 
proposed by Masao Iri, Sigeiti Moriguti, and Yoshimitsu Takasawa \cite{IMT70}.
The formula is known today as the ``IMT rule,''
which name was introduced in \cite{TM73}
and used in \cite{DR75}.

For an integral $I = \int_{0}\sp{1} f(x) {\rm d}x$
over $[0,1]$,
the IMT rule 
applies the trapezoidal rule to
the integral
$I = \int_{0}\sp{1} f(\phi (t)) \phi'(t) {\rm d}t$
resulting from the transformation 
\[
\phi(t) = \frac{1}{Q} \int_{0}\sp{t} 
   \exp  \left[ 
   - \left( \frac{1}{\tau} + \frac{1}{1-\tau} \right) \right]
 {\rm d}\tau, 
\]
where 
\[
Q =  \int_{0}\sp{1} 
   \exp  \left[ 
   - \left( \frac{1}{\tau} + \frac{1}{1-\tau} \right) \right]
 {\rm d}\tau
\]
is a normalizing constant
to render $\phi(1)=1$.

The transformed integrand 
$g(t) = f(\phi (t)) \phi'(t)$
has the property that all the derivatives 
$g^{(j)}(t)$ $(j=1,2,\ldots)$
vanish at $t=0,1$.
By the Euler--Maclaurin formula,
this indicates that the IMT rule should be highly accurate.
Indeed, it was shown in \cite{IMT70} via a complex analytic method
that the error of the IMT rule can be estimated roughly as 
$\exp ( -C \sqrt{N} )$,
which is much better than $N\sp{-4}$ of the Simpson rule, say, but  not as good as
$\exp ( -C N/ \log N )$
of the DE formula.
Variants of the IMT rule have been proposed for possible improvement
\cite{Mor78,MI82,Oou08IMT,Sug88},
but it turned out that an IMT-type rule,
transforming
$\int_{0}\sp{1} {\rm d}x$
to
$\int_{0}\sp{1} {\rm d}t$
rather than to
$\int_{-\infty}\sp{+\infty} {\rm d}t$,
cannot outperform the DE formula.

\section{DE-Sinc methods}

Changing variables is
also useful in the Sinc numerical methods.
The book by Stenger \cite{Ste93} in 1993 describes this methodology
to the full extent,
focusing on single exponential (SE) transformations like
$\phi(t) = \tanh (t/2)$.
Use of the double exponential transformation in the Sinc numerical methods
was initiated by Sugihara 
\cite{Sug97sinc,Sug03} around 2000,
with subsequent development mainly in Japan.
Such numerical methods are often called the DE-Sinc methods.
The subsequent results obtained in the first half of the 2000s are described 
in \cite{Mor05,MS01,SM04sinc}.

\subsection{Sinc approximation}

The Sinc approximation of a function $f(x)$ over $(-\infty,\infty)$
is given by 
\begin{equation} \label{Sincappr}
 f(x) \approx  \sum_{k=-N}\sp{N} f(kh) S(k,h)(x) ,
\end{equation}
where $S(k,h)(x)$ is the so-called Sinc function defined by 
\begin{equation}
  S(k,h)(x) = 
\frac{\sin [(\pi/h)(x - kh)]}{ (\pi/h)(x - kh) } 
\end{equation}
and the step size $h$ is chosen appropriately, depending on $N$.
The technique of variable transformation $x = \phi(t)$ 
is  also effective in this context.
By applying the formula \eqref{Sincappr} to $f(\phi(t))$ we obtain
\begin{equation}
 f(\phi(t)) \approx   \sum_{k=-N}\sp{N}  f(\phi(kh)) S(k,h)(t) ,
\end{equation}
or equivalently,
\begin{equation}
 f(x) \approx  \sum_{k=-N}\sp{N} f(\phi(kh)) S(k,h)(\phi\sp{-1}(x)) .
\end{equation}
To approximate $f(x)$ over $[0,1]$, for example, we choose
\begin{align}
 \phi(t) &= \frac{1}{2} \tanh \frac{t}{2} + \frac{1}{2} ,
\label{transSEsinc}
\\
 \phi(t) &= \frac{1}{2} \tanh \left( \frac{\pi}{2} \sinh t \right) + \frac{1}{2} ,
\label{transDEsinc}
\end{align}
etc.
The methods using
\eqref{transSEsinc} and
\eqref{transDEsinc} are often called the SE- and DE-Sinc approximations,
respectively.
When the approximation error is measured by
\begin{equation} \label{errorSincApprox}
\sup_{x \in [0,1]} 
\left|  f(x) - \sum_{k=-N}\sp{N} f(\phi(kh)) S(k,h)(\phi\sp{-1}(x)) \right| ,
\end{equation}
the error of the SE-Sinc approximation is roughly 
$\exp( -C \sqrt{N})$
and that of the DE-Sinc approximation is  
$\exp( -C N / \log N )$.

\begin{figure}
\centering
\includegraphics[width=0.6\textwidth]{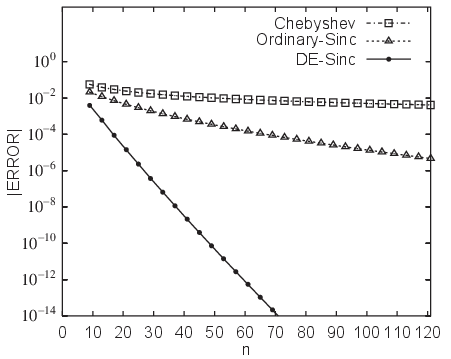}
\caption{
Errors in the Sinc approximations for $x\sp{1/2} (1-x)\sp{3/4}$
using \eqref{transSEsinc} and \eqref{transDEsinc} and the Chebyshev interpolation.  
``Ordinary-Sinc"
means the SE-Sinc approximation using \eqref{transSEsinc}.
The figure is taken from 
Sugihara and Matsuo \cite[Fig.~3]{SM04sinc}
with permission from Elsevier;
$n$ in the figure correspond to $N$ in \eqref{Sincappr}.
}
\label{FGdesincSM04}
\end{figure}

These approximation schemes
are compared in Fig.~\ref{FGdesincSM04}
(taken from \cite{SM04sinc})
for function $f(x)=x\sp{1/2} (1-x)\sp{3/4}$ over $[0,1]$.
The approximation error \eqref{errorSincApprox},
which is defined as the supremum of the discrepancy over $x \in [0,1]$,
is estimated by the maximum of 
the discrepancy at 2000 equally spaced points $x$ in $[0,1]$.
In Fig.~\ref{FGdesincSM04},  
``Ordinary-Sinc"
means the SE-Sinc approximation
using \eqref{transSEsinc},
and the polynomial interpolation 
with the Chebyshev nodes is included for comparison.

Detailed theoretical analyses on the DE-Sinc method
can be found in
Sugihara \cite{Sug03}
as well as 
Tanaka et al.~\cite{TSM09succDE} and Okayama et al.~\cite{OMS13,OTMS13}.
An optimization technique is used to improve 
the DE-Sinc method in Tanaka and Sugihara \cite{TS21const}.

\subsection{Application to other problems}

Once a function approximation scheme is at hand,
we can apply it to a variety of numerical problems.
Indeed this is also the case with the DE-Sinc approximation
as follows.

\begin{itemize}
\item
Indefinite integration by 
Muhammad and Mori \cite{MuM03},   
Tanaka et al.~\cite{TSM04indefint},
and Okayama and Tanaka~\cite{OT22}.

\item
Initial value problem of differential equations by 
Nurmuhammad et al.~\cite{NuMuM05ini}
and Okayama \cite{Oka18}.

\item
Boundary value problem of differential equations
by Sugihara \cite{Sug02},
followed by
Nurmuhammad et al.~\cite{NuMuM07galer,NuMuMS05bndry}
and Mori et al.~\cite{MNuMu09}.

\item
Volterra integral equation
by Muhammad et al.~\cite{MuNuMS05inteq}
and 
Okayama et al.~\cite{OMS15}.

\item
\hyphenation{Koba-yashi}
Fredholm integral equation by 
Kobayashi et al.~\cite{KOZ03}, 
Muhammad et al.~\cite{MuNuMS05inteq}, and
Okayama et al.~\cite{OMS11}.

\item
Computation of matrix functions (logarithm, fractional power)
by Tatsuoka et al.~\cite{TSMKZ21matfrac,TSMZ20matlog}. 

\end{itemize}


\noindent {\bf Acknowledgement}.
The authors are grateful to Ken'ichiro Tanaka 
and Tomoaki Okayama for their 
support in writing this article.
This work was supported by JSPS/MEXT KAKENHI 
23K21674 and 24K22290.




\begin{thebibliography}{99}

\bibitem{DR75}   
Davis, P.J.,   Rabinowitz, P.: 
Methods of Numerical Integration,
1st edn. Academic Press, NewYork (1975); 2nd edn.\ (1984)


\bibitem{IMT70}   
Iri, M.,   Moriguti, S.,  Takasawa, Y.:
On a certain quadrature formula (in Japanese).
{RIMS Kokyuroku} {\bf 91}, 82--118  (1970);
English translation in
{J. Comput. Appl. Math.}
{\bf 17}, 3--20  (1987) 



\bibitem{KOZ03}   
Kobayashi, K.,   Okamoto, H.,  Zhu, J.: 
Numerical computation of water and solitary waves by the double exponential transform.
{J. Comput. Appl. Math.}
{\bf 152}, 229--241 (2003)






\bibitem{Mor78}
Mori, M.:
An IMT-type double exponential formula for numerical integration.
{Publ. RIMS}
{\bf 14}, 713--729 (1978)






\bibitem{Mor05}
Mori, M.:
Discovery of the double exponential transformation and its developments.
{Publ. RIMS}
{\bf 41}, 897--935 (2005)



\bibitem{MNuMu09}
Mori, M.,   Nurmuhammad, A.,  Muhammad, M.:
DE-sinc method for second order singularly perturbed boundary value problems.
{Japan J. Indust. Appl. Math.}
{\bf 26}, 41--63  (2009)




\bibitem{MS01}
Mori, M.,  Sugihara, M.: 
The double exponential transformations in numerical analysis.
{J. Comput. Appl. Math.}
{\bf 127}, 287--296 (2001)



\bibitem{MuM03}   
Muhammad, M.,  Mori, M.:
Double exponential formulas for numerical indefinite integration.
{J. Comput. Appl. Math.}
{\bf 161}, 431--448  (2003)



\bibitem{MuNuMS05inteq}   
Muhammad, M.,  Nurmuhammad, A.,   Mori, M.,  Sugihara, M.:
Numerical solution of integral equations by means of 
the Sinc collocation method based on the double exponential transformation. 
{J. Comput. Appl. Math.}
{\bf 177}, 269--286  (2005)




\bibitem{MI82}   
Murota, K.,   Iri, M.:
Parameter tuning and repeated application of the IMT-type
transformation in numerical quadrature.
{Numer. Math.}
{\bf 38}, 347--363  (1982)

\bibitem{NuMuM05ini}    
Nurmuhammad, A.,   Muhammad, M.,  Mori, M.:
Numerical solution of initial value problems based on the double exponential transformation.
{Publ. RIMS}
{\bf 41}, 937--948  (2005)


\bibitem{NuMuM07galer}   
Nurmuhammad, A.,   Muhammad, M.,  Mori, M.:
Sinc-Galerkin method based on the DE transformation for the boundary value problem 
of fourth-order ODE.
{J. Comput. Appl. Math.}
{\bf 206}, 17--26  (2007)



\bibitem{NuMuMS05bndry}  
Nurmuhammad, A.,   Muhammad, M., Mori, M.,   Sugihara, M.:
Double exponential transformation in the Sinc-collocation method 
for a boundary value problem with fourth order ordinary differential equation.
{J. Comput. Appl. Math.}
{\bf 182}, 32--50  (2005)





\bibitem{OS05deYears}
Okamoto, H.,  Sugihara, M.,   eds.:
Thirty Years of the Double Exponential Transforms.
Special issue of 
{Publ. RIMS}
{\bf 41}, Issue 4 (2005)




\bibitem{Oka18}  
Okayama, T.:
Theoretical analysis of Sinc-collocation methods and 
Sinc-Nystr{\" o}m methods for systems of initial value problems.
{BIT Numer. Math.} 
{\bf 58}, 199--220  (2018)




\bibitem{OMS13}
Okayama, T.,   Matsuo, T.,  Sugihara, M.:
Error estimates with explicit constants for Sinc approximation, 
Sinc quadrature and Sinc indefinite integration.
{Numer. Math.} 
{\bf 124}, 361--394  (2013)


\bibitem{OMS11}    
Okayama, T.,   Matsuo, T., Sugihara, M.:
Improvement of a Sinc-collocation method for Fredholm integral equations of the second kind.
{BIT Numer. Math.} 
{\bf 51}, 339--366  (2011)



\bibitem{OMS15}  
Okayama, T.,   Matsuo, T.,   Sugihara, M.:
Theoretical analysis of Sinc-Nystr{\" o}m methods for Volterra integral equations.
{Math.\ Comput.}
{\bf 84}, 1189--1215  (2015)



\bibitem{OT22}   
Okayama, T.,    Tanaka, K.:
Yet another DE-Sinc indefinite integration formula.
{Dolomites Res. Notes Approx.}
{\bf 15}, 105--116  (2022)






\bibitem{OTMS13}   
Okayama, T.,   Tanaka, K.,  Matsuo, T.,  Sugihara, M.:
DE-Sinc methods have almost the same convergence property as SE-Sinc methods 
even for a family of functions fitting the SE-Sinc methods, 
Part I: definite integration and function approximation.
{Numer. Math.} 
{\bf 125}, 511--543   (2013)




\bibitem{Oou08IMT}   
Ooura, T.:
An IMT-type quadrature formula with the same asymptotic performance 
as the DE formula.
{J. Comput. Appl. Math.}
{\bf 213}, 232--239   (2008)



\bibitem{OM91osc}   
Ooura, T.,  Mori, M.:
The double exponential formula for oscillatory functions over the half infinite interval.
{J. Comput. Appl. Math.}
{\bf 38}, 353--360  (1991)



\bibitem{OM99four}   
Ooura, T.,  Mori, M.:
A robust double exponential formula for Fourier type integrals. 
{J. Comput. Appl. Math.}
{\bf 112}, 229--241   (1999)



\bibitem{Sch69} 
Schwartz, C.:
Numerical integration of analytic functions.
{J. Comput. Phys.} 
{\bf 4}, 19--29   (1969)





\bibitem{Ste73} 
Stenger, F.:
Integration formulas based on the trapezoidal formula.
{J. Inst. Math. Appl.} 
{\bf 12}, 103--114   (1973)


\bibitem{Ste78} 
Stenger, F.:
Optimal convergence of minimum norm approximations in  $H_p$.
{Numer. Math.}
{\bf 29}, 345--362  (1978)


\bibitem{Ste93}   
Stenger, F.:
Numerical Methods Based on Sinc and Analytic Functions.
Springer, New York (1993)



\bibitem{Sug86}   
Sugihara, M.: 
On optimality of the double exponential formulas (in Japanese).
{RIMS Kokyuroku} 
{\bf 585}, 150--175   (1986)




\bibitem{Sug88}   
Sugihara, M.: 
On optimality of the double exponential formulas, II (in Japanese).
{RIMS Kokyuroku} 
{\bf 648}, 20--38   (1988)





\bibitem{Sug97optDE}   
Sugihara, M.: 
Optimality of the double exponential formula---functional analysis approach.
{Numer. Math.}
{\bf 75}, 379--395   (1997)



\bibitem{Sug97sinc}   
Sugihara, M.: 
Sinc approximation using double exponential transformations (in Japanese).
{RIMS Kokyuroku} 
{\bf 990}, 125--134    (1997)




\bibitem{Sug02}  
Sugihara, M.: 
Double exponential transformation in the Sinc-collocation method
for two-point boundary value problems.
{J. Comput. Appl. Math.}
{\bf 149}, 239--250  (2002)


\bibitem{Sug03}
Sugihara, M.: 
Near optimality of the sinc approximation.
{Math.\ Comput.}
{\bf 72}, 767--786   (2003)




\bibitem{SM04sinc}
Sugihara, M.,  Matsuo, T.: 
Recent developments of the Sinc numerical methods.
{J. Comput. Appl. Math.}
{\bf 164/165}, 673--689   (2004)




\bibitem{TM73}
Takahasi, H.,  Mori, M.:
Quadrature formulas obtained by variable transformation.
{Numer. Math.}
{\bf 21}, 206--219    (1973)



\bibitem{TM73defirst}
Takahasi, H.,  Mori, M.:
Quadrature formulas obtained by variable transformation (2) (in Japanese).
{RIMS Kokyuroku} 
{\bf 172}, 88--104   (1973)








\bibitem{TM74}
Takahasi, H.,  Mori, M.:
Double exponential formulas for numerical integration.
{Publ. RIMS}
{\bf 9}, 721--741   (1974)



\bibitem{TO23bk}
Tanaka, K.,  Okayama, T.: 
Numerical Methods with Variable Transformations
(in Japanese).
Iwanami, Tokyo (2023)


\bibitem{TS21const}
Tanaka, K.,  Sugihara, M.: 
Construction of approximation formulas 
for analytic functions by mathematical optimization.
In: Baumann, G. (ed.) 
New Sinc Methods of Numerical Analysis,
pp.~341--368, Birkh{\" a}user, Basel (2021)  



\bibitem{TSM04indefint}
Tanaka, K.,    Sugihara, M.,  Murota, K.: 
Numerical indefinite integration by double exponential sinc method.
{Math.\ Comput.}
{\bf 74}, 655--679   (2004)




\bibitem{TSM09succDE}
Tanaka, K.,    Sugihara, M.,   Murota, K.: 
Function classes for successful DE-Sinc approximations.
{Math.\ Comput.}
{\bf 78}, 1553--1571    (2009)


\bibitem{TSMM09fnclass}
Tanaka, K.,    Sugihara, M.,   Murota, K.,   Mori, M.: 
Function classes for double exponential integration formulas. 
{Numer. Math.}
{\bf 111}, 631--655    (2009)


\bibitem{TSMKZ21matfrac}
Tatsuoka, F.,  Sogabe, T., Miyatake, Y.,  Kemmochi, T.,  Zhang, S.-L.: 
Computing the matrix fractional power based on the double exponential formula.
{Electron. Trans. Numer. Anal.} {\bf 54}, 558--580  (2021) 



\bibitem{TSMZ20matlog}
Tatsuoka, F.,  Sogabe, T., Miyatake, Y., Zhang, S.-L.: 
Algorithms for the computation of the matrix logarithm based 
on the double exponential formula.
{J. Comput. Appl. Math.} {\bf 373},  112396 (2020) 


\end{thebibliography}
\end{document}